\documentclass[11pt]{amsart}
\usepackage{amsmath}
\usepackage{amsfonts}
\usepackage{amssymb}
\usepackage{epsfig}

\newcommand{\nula}{\nu_\lambda}
\newcommand{\nuhat}{\widehat{\nu}}
\newcommand{\nulahat}{\widehat{\nu}_\lam}

\newcommand{\be}{\begin{eqnarray}}
\newcommand{\ee}{\end{eqnarray}}

\newcommand{\dist}{\mbox{\rm dist}}

\newcommand{\const}{\mbox{\rm const}}

\newcommand{\clos}{\mbox{\rm clos}}

\newcommand{\Int}{{\rm int}}
\newcommand{\es}{\emptyset}
\newcommand{\half}{\frac{1}{2}}
\newcommand{\shalf}{\textstyle{\half}}
\newcommand{\Leb}{{\mathcal L}}

\newcommand{\e}{{\varepsilon}}
\newcommand{\Th}{{\theta}}
\newcommand{\RE}{{\rm Re}}
\newcommand{\IM}{{\rm Im}}

\newcommand{\R}{{\mathbb R}}

\newcommand{\Z}{{\mathbb Z}}
\newcommand{\Nat}{{\mathbb N}}
\newcommand{\Compl}{{\mathbb C}}
\newcommand{\Comp}{{\mathbb C}}
\newcommand{\ov}{\overline}
\newcommand{\disc}{{\mathbb D}}
\newcommand{\Ak}{{\mathcal A}}

\newcommand{\ba}{{\bf a}}
\newcommand{\bx}{{\bf x}}

\newcommand{\lam}{\lambda }

\newcommand{\eps}{\varepsilon}

\newcommand{\Hau}{{\mathcal H}}

\newcommand{\Bk}{{\mathcal B}}

\newcommand{\Mk}{{\mathcal M}}
\newcommand{\Pk}{{\mathcal P}}
\newcommand{\Mtil}{\widetilde{\Mk}}

\newcommand{\QED}{\qed}

\newtheorem{theorem}{Theorem}[section]
\newtheorem{lemma}[theorem]{Lemma}

\newtheorem{prop}[theorem]{Proposition}

\theoremstyle{definition}
\newtheorem{defi}[theorem]{Definition}

\theoremstyle{remark}

\numberwithin{equation}{section}

\begin{document}

\title[]{On the ``Mandelbrot set'' for a pair of linear
maps and complex Bernoulli convolutions}

\author{Boris Solomyak}

\thanks{The authors were supported in part by NSF grant \#DMS 0099814 \\
\indent 2000 {\em Mathematics Subject Classification.} Primary 28A80;
Secondary 28A78, 37F45, 11R06, 26C10}

\author{Hui Xu}

\address{Department of Mathematics,
University of Washington, Seattle, WA 98195-4350 \newline
\indent {\tt solomyak@math.washington.edu} \newline
\indent {\tt hxu@math.washington.edu}} 

\begin{abstract}
We consider the family of self-similar sets $A_\lam$, attractors of the
iterated function system $\{\Compl;\ \lam z-1, \lam z+1\}$,  
depending on a parameter $\lam$ in the open unit disk.
First we study the set $\Mk$ of those $\lam$ for which $A_\lam$ is 
connected.  We show that a non-trivial portion of $\Mk$
near the imaginary axis is the closure of its interior (it is
conjectured that $\Mk\setminus \R$ is contained in the closure of its interior).
Next we turn to the sets $A_\lam$ themselves and natural measures
$\nula$ supported on them. These measures are the complex analogs of
much-studied infinite Bernoulli convolutions. 
Extending the results of Erd\H{o}s and Garsia, we demonstrate how 
certain classes of complex algebraic integers give rise to singular and
absolutely continuous measures $\nula$. Next we investigate
the Hausdorff dimension and measure 
of $A_\lam$, for Lebesgue-a.e.\ $\lam\in \Mk$, and 
obtain partial results on the absolute continuity of
$\nula$, for a.e.\ $\lam$ with $|\lam|>1/\sqrt{2}$.
\end{abstract}

\maketitle

\thispagestyle{empty}


\section{Introduction}
\label{sec-intro}

Consider a family of iterated function systems (IFS) in the complex
plane $\{\Compl;\ \lam z-1, \lam z+1\}$ depending on a parameter
$\lam \in \disc:= \{z\in \Compl:\ |z| < 1\}$. Let $A_\lam$ denote
the attractor of the IFS, that is, $A_\lam$ is the unique non-empty
compact set such that
\be \label{atra}
A_\lam = (-1+\lam A_\lam) \cup (1+\lam A_\lam).
\ee
These are among the most basic self-similar sets in the plane, and there
is considerable interest in understanding their topological and ``fractal''
properties. 

There is a fundamental dichotomy for attractors of IFS with two maps:
they are either connected or totally disconnected, depending on whether
the sets in the right-hand side of (\ref{atra}) have non-empty intersection.
Consider the set
$$
\Mk:= \{\lam\in \disc:\ A_\lam \mbox{\ is connected}\}\,.
$$
It was first introduced by Barnsley and Harrington \cite{BH}, who called it the
``Mandelbrot set for the pair of linear maps,'' by analogy with the
classical Mandelbrot set in complex dynamics. In fact, there are some
parallels with the quadratic family $z^2+c$, 
the attractors $A_\lam$ being analogs of Julia
sets, see \cite{bandt}.  The set $\Mk$ was studied in 
\cite{bou,IKR,IJK,solo,bandt}, but some
of the basic questions on the geometry and topology of $\Mk$ remain open.
See Figure 1, made by C. Bandt, which shows the part of $\Mk$ in
$\{z:\ \RE(z)>0,\,\IM(z)>0,\,|z| \le 1/\sqrt{2}\}$.
Computer pictures
suggest many interesting features, but it is challenging to prove
them rigorously. Our first result is a step toward proving that
$\Mk\setminus \R$ is contained in the closure of its interior, as conjectured by
Bandt \cite{bandt}.

\begin{figure}[ht]
\epsfig{figure=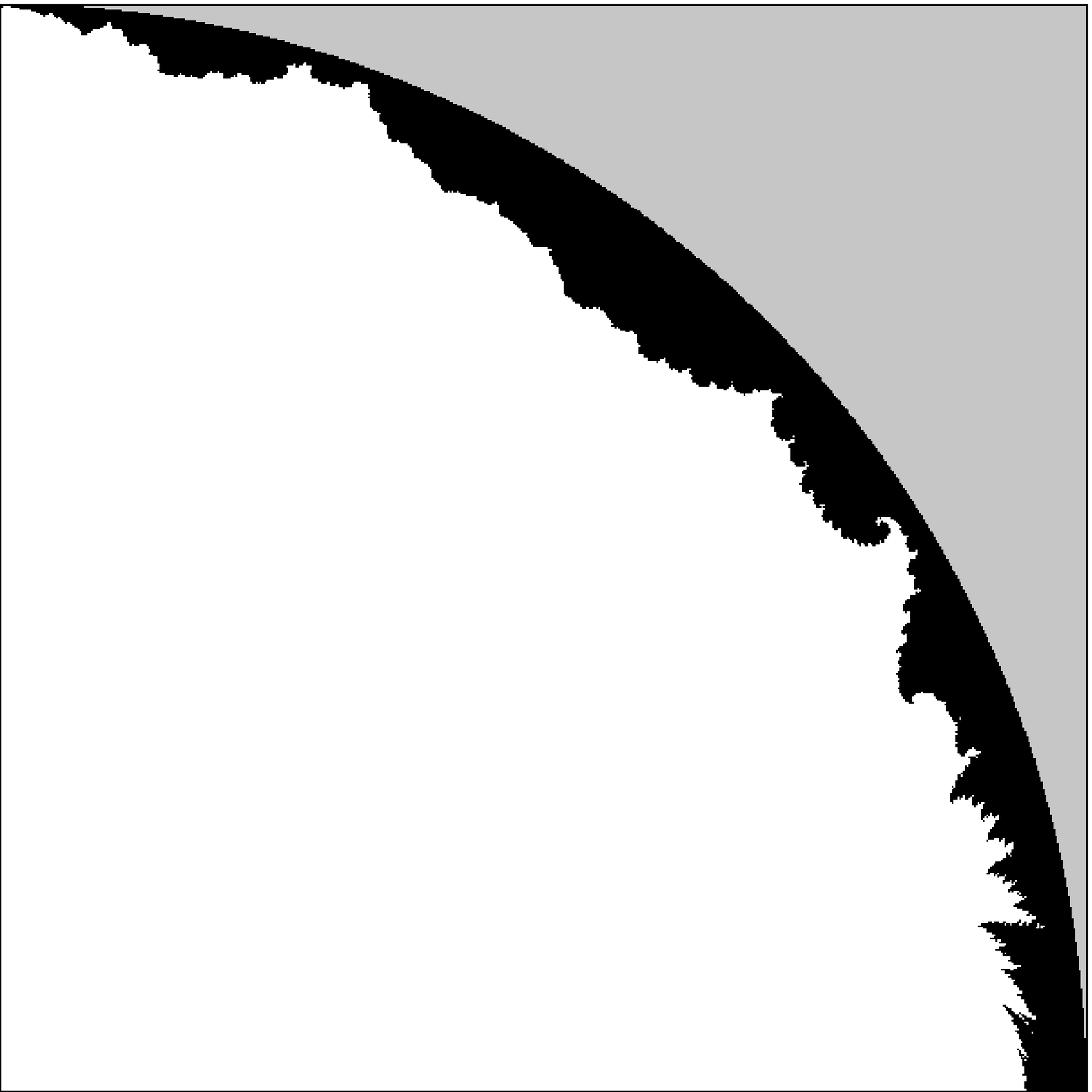,height=12cm}
\caption{The ``Mandelbrot set'' $\Mk$}
\end{figure}

Next we turn to the attractors $A_\lam$. The case
$\lam\not\in \Mk$ is completely understood, since then $A_\lam$ is
totally disconnected and has Hausdorff dimension equal to the
similarity dimension $s(\lam):= \frac{\log 2}{-\log|\lam|}$, with the
corresponding Hausdorff measure positive and finite. Thus, we focus our
attention on the attractors
$A_\lam$ for $\lam\in \Mk$. Then, at least for a ``typical'' 
$\lam$, the set $A_\lam$ has an ``overlap.'' 
Along with the set $A_\lam$ we consider the {\em natural measure}
$\nula$ supported on it, obtained by assigning
equal probabilities $\half$ to each of the maps in the IFS. These measures
are the complex analogs of {\em infinite Bernoulli convolutions}, much
studied since the 1930's, see \cite{sixty} for a survey.
We extend classical results of Erd\H{o}s \cite{erdos1} and
Garsia \cite{garsia1} to obtain two 
classes of complex algebraic integers that give rise to singular and
absolutely continuous measures $\nula$ respectively.

Our final theme is the so-called ``typical $\lam$'' or ``almost sure results.''
Since the sets $A_\lam$ and the measures $\nula$, for $\lam \in \Mk$, are
so difficult to study, we attempt to establish dimension
formulae, etc., for {\em almost every} $\lam$ with respect to the area
measure $\Leb_2$. This line of research has been very active in the last
10 years or so, see, e.g., \cite{polli,solerd,pesolo2,peschl}. 
They key to obtain results
of this type is a certain {\em transversality condition}, which in our case
means the absence of double zeros for power series with $\{0,\pm 1\}$
coefficients.  The extension to the
complex plane was started in \cite{solo}, but it was hampered by the
lack of effective estimates of double zeros. Recently, we became aware
of such estimates in \cite{bbbp}, which yield definitive
``typical $\lam$'' statements for $\lam \in \Mk$, with $|\lam|< 1/\sqrt{2}$
(when $\nula$ is {\em a priori} singular),
and partial results for $|\lam|> 1/\sqrt{2}$ (when $\nula$ is expected to
be absolutely continuous for a typical $\lam$). 

\medskip

To conclude the introduction, we note that the sets $A_\lam$ and the
measures $\nula$ arise as ``building blocks'' when analyzing attractors
$E(T,\ba)$
of IFS $\{T\bx, T\bx + \ba\}$ where $T$ is a linear transformation
in $\R^d$ and $\ba$ is a generic non-zero vector. For instance, if $T$ has
a complex eigenvalue $\lam$, then the projection of $E(T,\ba)$
onto the $T$-invariant real plane corresponding to $\lam$ is an affine
copy of $A_\lam$. This served as an additional motivation for our study.


\section{Statement of results}

\subsection{Structure of $\Mk$.} We begin with some basic facts about 
$A_\lam$ and $\Mk$; 
they are not difficult and may be found in \cite{BH,barnsley}.
Then we briefly mention more recent results; see \cite{bandt} for more 
information.
We have
\be\label{eq-ss2}
A_\lam = \Big\{\sum_{n=0}^\infty a_n \lam^n:\ a_n \in
\{-1,1\}\Big\}\,.
\ee
Since $\lam \in \Mk$ if and only if $(\lam A_\lam + 1) \cap (\lam A_\lam-1)
\ne \es$, it follows that $\Mk$ is the set of zeros of $\{0,\pm 1\}$ power series:
\be\label{eq-mand1}\ \ \ \ \ \ \ 
\Mk = \Big\{\lam \in \disc:\ \exists\,\{a_k\}_{k=1}^\infty,
\ a_k \in \{-1,0,1\},\
1 + \sum_{k=1}^\infty a_k \lam^k = 0 \Big\}\,.
\ee
From (\ref{eq-mand1}) it follows
that $\Mk$ is relatively closed in the unit disc $\disc$.
The following provides easy estimates for $\Mk$ in terms of $|\lam|$:
\be\label{eq-inclu}
\{\lam\in \disc:\ |\lam| \ge 2^{-1/2}\} \subset \Mk \subset
\{\lam\in \disc:\ |\lam| \ge 1/2\}.
\ee

In their very interesting paper, Odlyzko and Poonen \cite{OP}
investigated the set of zeros of power series with coefficients 0,1.
Although \cite{OP}
does not have an immediate application to $\Mk$ and
the family $\{A_\lam\}$, it turned out to be very useful in this area.
In particular, using the ideas of \cite{OP},
Bousch \cite{bou} proved that $\Mk$ is connected and
locally connected (in contrast with the classical Mandelbrot set, for which
local connectedness is a famous open problem).
Interesting new results were recently obtained by 
Bandt \cite{bandt}.  In particular, he gave a rigorous 
computer-assisted proof that
$\Mk$ has ``holes'' (i.e., $\disc \setminus \Mk$ has more than one component).

A peculiar feature of the set $\Mk$ is the ``antenna'' (or ``spike'')
on the positive real axis, from $0.5$ to about $.67$, see \cite{BH,solo,bandt}.
More precisely, there is a line segment $J=[.5,\alpha]$, with
$\alpha \approx .67$, such that 
$
J \not \subset \clos(\Mk \setminus \R).
$
(Of course, there is a symmetric ``antenna,'' on the negative real axis.
The set $\Mk$ is clearly symmetric with respect to both axes, so we will
always confine ourselves to the first quarter of the plane.)
By (\ref{eq-inclu}), the ``interesting part'' of the set 
$\Mk$ lies in the disc $\{\lam \in \disc:\ |\lam| \le 2^{-1/2}\}$.
Bandt \cite{bandt} conjectured that the set $\Mk\setminus \R$ 
is contained in the closure of $\Int(\Mk)$. In the next theorem we prove
a partial result in this direction.

\medskip

\noindent {\bf Notation.}
Denote by $B_r(z_0)$ the open disc of radius $r$ centered at $z_0$ and write
$B_r:= B_r(0)$.

\begin{theorem} \label{th-main}
Let $H= \{\lam\in B_{1/\sqrt{2}}:\ \RE(\lam)>0, \ \IM(\lam)>0\}
\setminus B_{2/3}(1/3)$; then
$$
\Mk \cap \Int(H) \subset \clos(\Int(\Mk)).
$$
\end{theorem}

The set $H$ is shown in Figure 2. One can check that
$$
H = \{\lam\in \disc:\ 1/3 \le |\lam|^2 \le 1/2,\ 
0 \le \RE(\lam) \le (3|\lam|^2 -1)/2 \}.
$$
The right-most point in $H$ is
$\frac{1}{4} + i \frac{\sqrt{7}}{4}$ (incidentally, the corresponding set
$A_\lam$ is known as the ``tame twindragon'').
Figure 3 (made by C. Bandt) shows 
the set
$\Mk \cap \{\lam:\ 0 \le \RE(\lam) \le 1/4\}$, 
which is slightly larger than $\Mk \cap H$.

\begin{figure}[ht]
\epsfig{figure=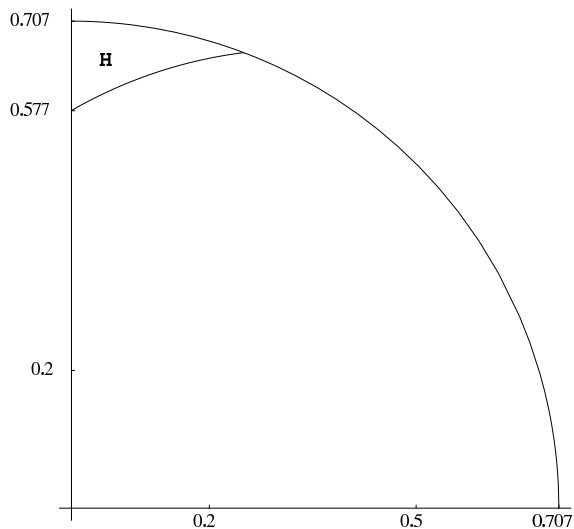,height=10cm}
\caption{The set $H$}
\end{figure}

\begin{figure}[ht]
\epsfig{figure=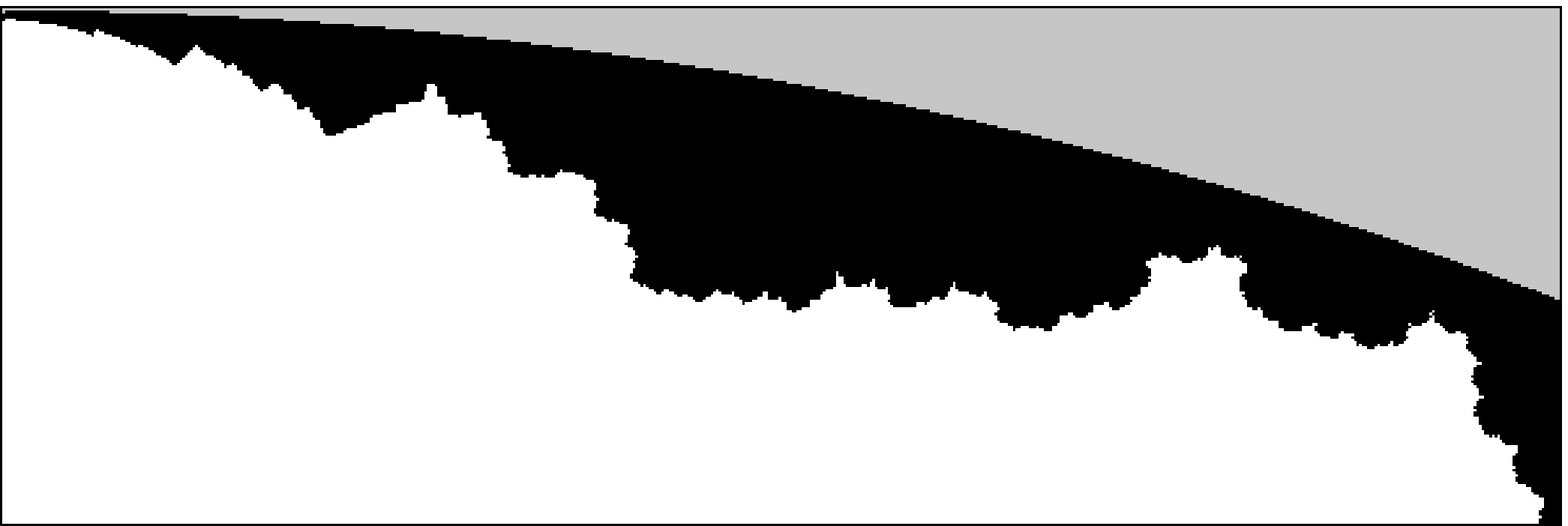,height=4.2cm}
\caption{A part of the set $\Mk$}
\end{figure}

Our method yields many concrete examples of interior
points in $\Mk \cap B_{1/\sqrt{2}}$. 
Let $\Mk_0$ be the set of zeros of {\em polynomials} with coefficients
in $\{0,\pm 1\}$, that is,
$$
\Mk_0:= \Big\{\lam \in \disc:\ \exists\,\{a_k\}_{k=1}^n,\ a_k \in \{-1,0,1\},\
1 + \sum_{k=1}^n a_k \lam^k = 0 \Big\}\,.
$$
We prove that all points in $\Mk_0 \cap \Int(H)$ lie in the
interior of $\Mk$. For example, $\lam_0 \approx .141964 + .677696 i$,
a zero of the polynomial $1+z^2 + z^3 - z^4 - z^5 + z^7$,
has this property, and we check that its neighborhood of radius $2\cdot
10^{-3}$ lies in $\Mk$.

We should note that our method is rather similar to that of Indlekofer, J\'arai and K\'atai,
who gave a computer-assisted proof that $\Mk$ contains
$e^{-2\pi i/5}/\sqrt{2}$ in its interior \cite[p.\,423]{IJK}, but there are
significant differences as well, see the discussion at the end of Section 3.

\subsection{Complex Bernoulli convolutions.}
Let $\nula$ be the distribution of the random series
$\sum_0^\infty\pm\lambda^n$ where the signs
are chosen independently with probabilities $\{\half,\half\}$.
It is  the infinite convolution product of $\half(\delta_{-\lam^n} +
\delta_{\lam^n})$, a probability measure on the plane. 
Alternatively, 
$\nula$ may be defined as the unique probability measure such that
$$
\nula = \shalf (\nula \circ S_1^{-1} + \nula \circ S_2^{-1}),
$$
where $S_1(z) = \lam z + 1$ and $S_2(z) = \lam z-1$, see \cite{hutch}.
Thus, $\nula$ is a
self-similar measure, and by the Law of Pure Type it is
either singular or absolutely continuous with respect to the Lebesgue
measure $\Leb_2$ (see \cite{sixty} for a short proof).
Observe that $A_\lam$ is precisely the compact support of $\nula$.
It is obvious that $\Leb_2(A_\lam) = 0$  for $|\lam| < 2^{-1/2}$,
so $\nula$ is singular. The question {\em ``for which $\lam$, with
$|\lam| \ge 2^{-1/2}$, is $\nula$ absolutely continuous ?''} is the complex
analog of a well-known problem, studied since the 1930's and still not
completely solved (see \cite{sixty} for a survey).

\begin{defi} \label{def-pisot}
An algebraic integer $\alpha>1$ is a {\em Pisot number} (or PV-number),
if all its Galois conjugates (i.e., other roots
of the minimal polynomial) are less than one in modulus.

A non-real algebraic integer $\theta$, with $|\Th|>1$,
is called a {\em complex Pisot number}
if all its Galois conjugates, except $\ov{\theta}$, are less than one in
modulus. 
\end{defi}

Below singular/absolutely continuous (a.c.)
is always understood with respect to the
planar Lebesgue measure $\Leb_2$, unless stated otherwise.

The following theorem extends the result of Erd\H{o}s \cite{erdos1},
who proved that $\nula$ is singular with respect to $\Leb_1$ for all
real $\lam\in (\half,1)$ such that $1/\lam$ is a Pisot number.
(It is an open problem whether this condition is also necessary, i.e., 
whether $\nula$ is a.c.\ with respect to $\Leb_1$ for all
real $\lam\in (\half,1)$ other than reciprocals of Pisot numbers.)

\begin{theorem} \label{thm-pisot}
If $\theta$ is a complex Pisot number and $1 < |\Th| < \sqrt{2}$,
then $\nula$ is singular for $\lam = 1/\Th$. 
\end{theorem}

\noindent
{\bf Remarks.} 1. If $|\Th|>\sqrt{2}$, then $\nula$, for $\lam = 1/\Th$,
is singular for the trivial reason that $\Leb_2(A_\lam)=0$. 

2. It is well-known (Siegel \cite{siegel}, see also \cite{pisot}) 
that the smallest real
Pisot number is $\alpha_0 \approx 1.3247$, the positive
zero of $z^3-z -1$. Chamfy \cite{chamfy} found the smallest in modulus complex
Pisot numbers: their modulus is $\sqrt{\alpha_0} \approx 1.1509$,
with either $z^3-z^2+1$ or $z^6-z^2+1$ as a minimal polymial.
Garth \cite{garth} found
a list of all sufficiently
small complex Pisot numbers. The reciprocals of ten smallest
complex Pisot numbers in the first quarter of the plane are shown in Figure 4.

\begin{figure}[ht]
\epsfig{figure=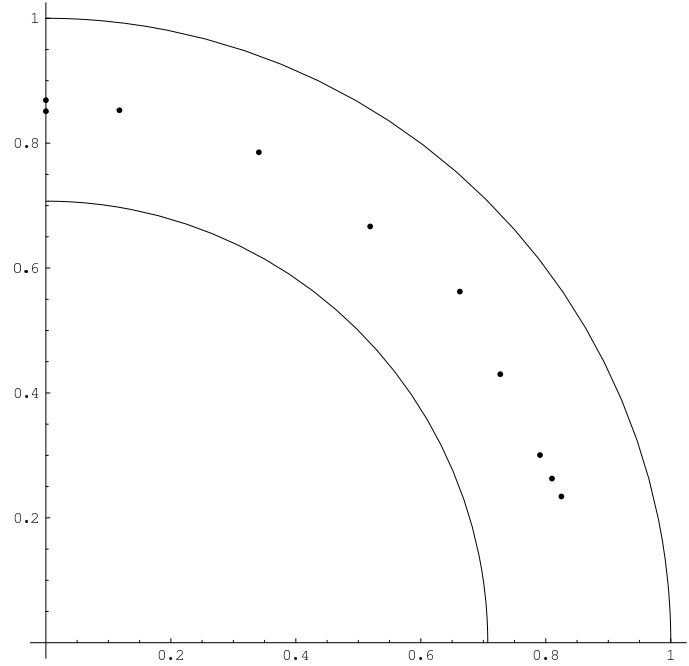,height=12cm}
\caption{Complex Pisot numbers}
\end{figure}

3. Let $S_1, S_2$ denote the sets of real and complex Pisot numbers 
respectively. There are no quadratic numbers in $S_2 \cap B_{\sqrt{2}}$.
Among the cubic numbers in $S_2 \cap B_{\sqrt{2}}$
are zeros of $z^3-z^2+z+1, z^3+z+1$, and $z^3-z^2+1$.
The set $S_2 \cap B_{\sqrt{2}}$ also includes all numbers
$\pm i \sqrt{\alpha}$ where $\alpha\in (1,2)$ is a real Pisot number.

4. The set $S_2 \cup \pm S_1$ is closed, see \cite[9.2]{pisot}.
 
\begin{defi}
We say that an algebraic integer $\theta$, with $|\theta|>1$,
is a {\em Garsia number}
if all its Galois conjugates are greater than one in modulus.
A {\em complex Garsia number} will mean a non-real Garsia number.
\end{defi}

The following theorem extends the result of Garsia \cite{garsia1},
who proved that $\nula$ is a.c.\ with respect to $\Leb_1$, with a bounded
density, for all
real $\lam\in (\half,1)$ such that $1/\lam$ is a Garsia number, whose
minimal polynomial has constant term $\pm 2$.

\begin{theorem} \label{th-garsia}
If $\theta$ is a complex Garsia number and the minimal polynomial
for $\Th$ has constant term $\pm 2$, then $\nula$, for $\lam = 1/\Th$,
is a.c.\ with a bounded density.
\end{theorem} 

\noindent {\bf Remarks.} 1. The conditions on $\Th$ imply that
$1  < |\Th| \le \sqrt{2}$, since $\ov{\Th}$ is always a conjugate and
the product of all zeros equals $\pm 2$. The equality
$|\Th|=\sqrt{2}$ implies that $\Th$ is quadratic.
There are exactly
ten such complex Garsia numbers: $\pm i\sqrt{2},\ \pm 1 \pm i,\
\pm \half \pm i\frac{\sqrt{7}}{2}$, which
yield three essentially different measures
$\nula$.
The compact supports $A_\lam$ of these measures are
well-known as ``reptiles'' that tile the plane periodically: 
$A_\lam$ is a rectangle
for $\lam = i\cdot 2^{-1/2}$, the ``twindragon'' for $\lam = \half+i\half$, and
the ``tame twindragon'' for $\lam = \frac{1}{4} + i \frac{\sqrt{7}}{4}$. 
The measure $\nula$ in each of these cases is just
the normalized Lebesgue measure restricted to $A_\lam$.

2. As in \cite{garsia1}, it is easy to see that if $q$ is
a polynomial with integer coefficients, having the leading coefficient one
and the constant term $\pm 2$, with all roots satisfying 
$|\Th_j| \in [1,2]$, then
every non-real root $\Th$ of $q$ is a complex Garsia number. (Indeed,
consider the minimal polynomial $p$ of $\Th$, of degree $\ell$.
It cannot have roots
of modulus one, since otherwise $p(z)$ and $z^\ell p(1/z)$ would have a common
zero, hence would have the same set of zeros, which is a contradiction.)
Thus, all non-real roots of the polynomials given in \cite[1.8]{garsia1}
are complex Garsia numbers. We get a countable set of such numbers
from the polynomials $x^{m+n} - x^n -2$, with $\max\{m,n\} \ge 2$.
It is easy to see that such a polynomial has all roots $|\Th_j| \in [1,2]$ and
at most two roots are real. 

3. Garsia numbers whose minimal polynomial has constant term $\pm m$, for
$m>2$,
occur in a similar context 
as well---their reciprocals $\lam$ have the property that the
random series $\sum_{j=0}^\infty a_j \lam^j$ has a.c.\ distribution,
where $a_j$ are random and independent and take values in
$\{0,1,\ldots,m-1\}$ (or any other set of $m$ equidistant digits)
with equal probabilities.
(If $\lam$ is real, then the distribution is
a.c.\ with respect to $\Leb_1$; otherwise it is a.c.\ 
with respect to $\Leb_2$.)

\subsection{The sets $A_\lam$ for $|\lam|\ge 1/\sqrt{2}$.}
If the measure $\nula$ is absolutely continuous, then, of course, 
its support
$A_\lam$ has positive area, and if the density of $\nula$ is continuous,
then $A_\lam$ has non-empty interior. However, the converse is false.
In particular, for all
$\lam \in [i\cdot 2^{-1/2},i)$ the set $A_\lam$ is a rectangle, and so 
obviously has non-empty interior, whereas
infinitely many $\lam$ in this segment on the imaginary axis are reciprocals
of complex Pisot numbers and so $\nula$ is singular.
It is plausible that $A_\lam$ has non-empty interior for {\em all}
non-real $\lam$ with $1/\sqrt{2} < |\lam| < 1$. Some elementary 
results in this direction are contained in the following proposition.
Part (ii) is a straightfoward extension of \cite[p.\,424]{IJK}.

\begin{prop} \label{prop-elem}

{\bf (i)} $A_\lam$ has non-empty interior for all $\lam$ in
\be
\Omega& := & \{\lam \in \disc:\ \RE(\lam)\ge 0,\ \IM(\lam)\ge 0\} \setminus
B_{\sqrt{3}/2}(1/2) \nonumber \\
      & = & \{\lam \in \disc:\ 0 \le \RE(\lam) \le |\lam|^2-(1/2)\}.
\nonumber \ee

{\bf (ii)} If $\lam = re^{\pi i m/n}$, with
$r \in [2^{-1/n},1)$, where $GCD(m,n)=1$, then
$A_\lam$ is a $(2n)$-gon having  all angles equal
to $\frac{\pi(n-1)}{n}$, so it has non-empty interior.
If $\lam = re^{2\pi i m/(2n+1)}$, with
$r \in [2^{-1/(2n+1)},1)$, where $GCD(m,2n+1)=1$, then
$A_\lam$ is a $(4n+2)$-gon having  all angles equal
to $\frac{\pi(2n)}{2n+1}$, so it has non-empty interior.
\end{prop}

For additional results on $A_\lam$ with non-empty interior, see
Theorem~\ref{th-meas}(ii) below.

\smallskip

For completeness, we mention what happens if $|\lam|=  1/\sqrt{2}$,
although we do not have anything new in this case. Then the similarity
dimension of $A_\lam$ equals 2, so $\Leb_2(A_\lam)>0$ if and only if
$A_\lam$ has non-empty interior \cite{schief}.
It is conjectured that this happens only
for quadratic $\lam$ in Remark 1 to Theorem~\ref{th-garsia}; 
this seems to be a difficult
problem. Some progress on a related conjecture was made in \cite{IJK}.

\medskip

\subsection{Results for a ``typical'' $\lam$.}
By ``typical'' we mean almost every (a.e.) with respect to the Lebesgue
measure $\Leb_2$. We do not review the history of this development, but refer
the reader to \cite{polli,solerd,pesolo2,solo,peschl,sixty}. A key
technical element needed for these results is a ``transversality condition.''
Let 
$$
\Bk_{\{0,\pm 1\}} = \Big\{1 + \sum_{k=1}^\infty a_k z^k:\ a_k \in \{-1,0,1\}
\Big\}\,.
$$
Thus, $\Mk = \{\lam\in \disc:\ \exists\,f\in \Bk_{\{0,\pm 1\}},\
f(\lam)=0\}$. Consider the set of {\em double zeros}
$$
\Mtil = \{\lam\in \disc:\ \exists\,f\in \Bk_{\{0,\pm 1\}},\
f(\lam)=f'(\lam)=0\}.
$$
We say that $\disc \setminus \Mtil$ is the {\em region of 
transversality} for $\Bk_{\{0,\pm 1\}}$. 
The following result was established in
1998 (and a preprint was available since 1995), although
it was not explicitly stated in this form.

\begin{prop} \label{prop-trans} {\bf \cite{bbbp}}
A power series of the form $1 + \sum_{n=1}^\infty a_n z^n$, with
$a_n\in [-1,1]$, cannot have a non-real double zero of modulus less than
$2\cdot 5^{-5/8} \approx .73143$.
\end{prop}

This yields the desired transversality on $\Mk \cap B_{1/\sqrt{2}}$
and in a thin annulus outside $B_{1/\sqrt{2}}$. Below we collect the
consequences for $A_\lam$ and $\nula$. 
Let
$$
s(\lam) = \frac{\log 2}{-\log|\lam|}
$$
denote the similarity dimension of the set $A_\lam$.

\begin{theorem} \label{th-dim}
We have $\dim_H(A_\lam) = s(\lam)$
for a.e.\ $\lam \in \Mk \cap  B_{1/\sqrt{2}}$.
Moreover, for any $0<r<R<1/\sqrt{2}$,
$$
\dim_H \{\lam \in \Mk:\ r < |\lam| < R,\ \dim_H A_\lam < s(\lam)\} \le
\frac{\log 2}{-\log R}<2.
$$
\end{theorem}

This follows from \cite{solo} and
\cite[Thm 8.2]{peschl}, in view of Proposition~\ref{prop-trans}.

\medskip

\noindent {\bf Remarks.} 1. Since $A_\lam$ is self-similar,
its Hausdorff, Minkowski, and packing dimensions are all equal \cite{falc3}.

2. It is well-known that $\dim_H(A_\lam)\le s(\lam)$ for all $\lam$,
see, e.g., \cite[9.6]{falcbook}.

3. It is proved in \cite[Prop 2.3]{solo} that
$\dim_H(A_\lam)< s(\lam)$ for a dense set of $\lam \in \Mk \cap
B_{1/\sqrt{2}}$. In fact, this holds for all $\lam \in \Mk_0 \cap
B_{1/\sqrt{2}}$. (Recall that $\Mk_0$ is the set of zeros of {\em polynomials}
with $\{0,\pm 1\}$ coefficients.)
It is a challenging open problem whether there exist $\lam \not\in \Mk_0
\cap B_{1/\sqrt{2}}$ such that $\dim_H(A_\lam)< s(\lam)$.

\medskip

Denote by $\Hau^s$ and $\Pk^s$ the $s$-dimensional Hausdorff and packing
measure respectively (see, e.g., \cite{falcbook} for definitions).
The next theorem says that for a typical $\lam \in \Mk\cap
B_{1/\sqrt{2}}$, the $s(\lam)$-dimensional Hausdorff measure of $A_\lam$
is zero,
whereas the packing measure $\Pk^{s(\lam)}(A_\lam)$
is positive and finite. 

\begin{theorem} \label{th-zero}

{\bf (i)} We have $\Hau^{s(\lam)} (A_\lam) = 0$ 
for a.e.\ $\lam \in \Mk \cap B_{1/\sqrt{2}}$.

{\bf (ii)} We have $\Pk^{s(\lam)} (A_\lam)< \infty$ for all $\lam$ and
$\Pk^{s(\lam)}(A_\lam)>0$ for a.e.\ $\lam \in \Mk \cap B_{1/\sqrt{2}}$.
Moreover, for any $0<r<R<1/\sqrt{2}$,
$$
\dim_H \{\lam \in \Mk:\ r < |\lam| < R,\ \Pk^{s(\lam)}(A_\lam)=0\} \le
\frac{\log 2}{-\log R}<2.
$$

{\bf (iii)} For a.e.\ $\lam \in \Mk \cap B_{1/\sqrt{2}}$ the measure
$\nula$ is equivalent to $\Pk^{s(\lam)}|_{A_\lam}$.
\end{theorem}

This is proved using the methods of \cite{pesiso}; again, transversality
(Proposition~\ref{prop-trans}) is essential.  The needed
modifications are fairly straightforward; they are left to the reader.

\medskip

\noindent {\bf Remarks.} 1. Schief \cite{schief} proved that
a self-similar set has zero Hausdorff measure in its similarity
dimension if and only if the {\em open set condition} fails. (In our
case the validity of the open set condition means, by definition, that there
is a non-empty open set $U\subset \Comp$ such that its images
$\lam U\pm 1$ lie in $U$ and do not intersect.) It is sometimes
considered that ``overlapping'' is synonymous with the failure
of the open set condition.

2. There is a ``topological'' version of the Hausdorff dimension and
Hausdorff measure result (but not of the packing measure result).
In fact, there is a dense $G_\delta$ subset
$\Phi \subset \Mk \cap B_{1/\sqrt{2}}$  such that $\dim_H(A_\lam)=
s(\lam)$, with $\Hau^{s(\lam)} (A_\lam) = 0$, for all $\lam \in
\Phi$. The dimension formula follows from the lower semi-continuity
of Hausdorff dimension as a function of parameter \cite{JV}, as in
\cite[Th.2.3]{SiSo}. Zero Hausdorff measure follows from \cite{pesiso}.

\medskip

Next we turn to the case $|\lam|>1/\sqrt{2}$ again.
Then the similarity dimension of $A_\lam$ exceeds 2. It is conjectured that
for a.e.\ such $\lam$ the measure $\nula$ is absolutely continuous and
the set $A_\lam$ has positive area. The
following theorem is a partial result in this direction.

\begin{theorem} \label{th-meas}
{\bf (i)} The measure $\nula$ is a.c.\ with a density in $L^2(\R^2)$, 
hence $\Leb_2(A_\lam)>0$, for a.e.\
$\lam \in \{z\in \Comp:\ 2^{-1/2} \le |z| \le 2\cdot 5^{-5/8}\}$.
(Note that $2^{-1/2} \approx .7071067,\ 2\cdot 5^{-5/8} \approx .7314316$.)
Moreover, for any $2^{-1/2}< r < R < 2\cdot 5^{-5/8}$,
$$
\dim_H \{\lam \in \Compl:\ r < |\lam| < R,\ d\nula/dx \not\in L^2(\R^2)\}
\le 4 - \frac{\log 2}{-\log r}\,.
$$

{\bf (ii)} The measure $\nula$ is a.c.\ with a continuous density, hence
$A_\lam$ has non-empty interior, for a.e.\ $\lam$ such that
$$
|\lam| \in \bigcup_{k=2}^\infty (2^{-1/(2k)}, (2\cdot 5^{-5/8})^{1/k}),
$$
in particular, for a.e.\ $\lam\in \disc$ with
$|\lam| > 2^{-1/20} \approx .9659363$.
\end{theorem}

Part (i) follows from \cite{solo} and
\cite[Thm 8.2]{peschl}, in view of Proposition~\ref{prop-trans}.
Part (ii) easily follows from (i);
we indicate the proof in Section 5.

\medskip

\noindent {\bf Remark.} Although with the methods of \cite{pesolo1,pesolo2}
we can somewhat increase the region where the
statements in the last theorem hold, 
unfortunately, we are not able to ``cover'' the whole
annulus $1/\sqrt{2} < |\lam| < 1$.

%
%


\section{Interior points of $\Mk$}

In this section we prove Theorem~\ref{th-main}.
The proof is based on several lemmas. They are all quite simple, but one of them
requires some calculations and its proof is postponed to the end of the
section. (We emphasize, however, that our proof {\em does not} rely
on computer.) Let
\be\label{eq-ss1}
A_\lam\{-1,0,1\} = \Big\{\sum_{n=0}^\infty a_n \lam^n:\ a_n \in
\{-1,0,1\}\Big\}\,.
\ee
Observe that $A_\lam\{-1,0,1\}$ is the attractor of the IFS
$\{\Compl;\ \lam z-1, \lam z, \lam z + 1\}$.
The following lemma is standard, see e.g. \cite[Lemma 7]{IJK}.

\begin{lemma} \label{lem-attr1}
If $F \subset \Compl$ is compact, $\lam \in \Compl$, and
\be \label{eq-attr1}
F \subset \lam F \cup (\lam F-1) \cup (\lam F + 1),
\ee
then $F \subset A_\lam\{-1,0,1\}$.
\end{lemma}

\begin{lemma} \label{lem-inclu}
If $F \subset A_\lam\{-1,0,1\}$ and there exists a finite sequence
$\{a_k\}_{k=1}^n$, $a_k\in \{-1,0,1\}$, such that
\be \label{eq-inclu2}
1 + \sum_{k=1}^n a_k \lam^k \in \lam^{n+1} F,
\ee
then $\lam \in \Mk$.
\end{lemma}

{\em Proof.} By (\ref{eq-ss1}), every point in $\lam^{n+1} F \subset
\lam^{n+1} A_\lam\{-1,0,1\}$ can be written as $\sum_{k=n+1}^\infty b_k \lam^k$,
with coefficients $b_k$ 
in $\{-1,0,1\}$. Thus, $\lam \in \Mk$ by (\ref{eq-mand1}).
\QED

\medskip

Let $R_{a,b}\subset \Compl$ denote the rectangle centered at the origin,
with the vertices at $\pm a \pm ib$. 

\begin{lemma} \label{lem-comput} 
For any $\lam\in H$, there
exist $a\ge 1$ and $b\ge 0.5$ such that 
$$
R_{a,b} \subset \lam R_{a,b} \cup (\lam R_{a,b} -1) \cup
(\lam R_{a,b} +1).
$$
\end{lemma}

The proof of this lemma is given at the end of the section.
Roughly speaking, the reason it works (at least in some
region) is that for $\lam\in [i 3^{-1/2},i)$, the attractor
$A_\lam\{-1,0,1\}$ is {\em exactly} a rectangle $R_{a,b}$ for appropriate
$a$ and $b$.
Recall that $\Mk_0$ denotes the set of zeros of 
{\em polynomials} with coefficients
in $\{0,\pm 1\}$. The following lemma is standard.

\begin{lemma} \label{lem-rouche}
$\Mk = \clos(\Mk_0)\cap \disc$.
\end{lemma}

{\em Proof.} Let $\lam\in \Mk$. We have 
$1 + \sum_{k=1}^\infty a_k \lam^k = 0$ for some $a_k \in \{-1,0,1\}$.
One can apply Rouch\'e's Theorem to show that for any $\e>0$ there 
exists $n$ such that the polynomial $1 + \sum_{k=1}^n a_k z^k$ has
a zero in $B_\e(\lam)$. \QED

\medskip

{\em Proof  of Theorem~\ref{th-main}.} 
It follows from Lemma~\ref{lem-attr1} and Lemma~\ref{lem-comput} that
$$F:=B_{0.5}(0) \subset A_\lam\{-1,0,1\}$$ for all $\lam \in H$.
In view of Lemma~\ref{lem-rouche}, it is enough to show that any
$\lam_0\in \Int(H)\cap \Mk_0$ is an interior point of $\Mk$.
Since $\lam_0\in \Mk_0$, there exist $a_1,\ldots,a_n\in \{-1,0,\}$ such that
$$
p(\lam_0)=0, \ \ \ \mbox{where}\ \ p(z) =  1 + \sum_{k=1}^n a_k z^k.
$$
We can find $\delta>0$ such that $B_\delta(\lam_0) \subset H$ and,
by the continuity of $p$, 
$$
|\lam-\lam_0| < \delta \ \Rightarrow\ |p(\lam)|=|p(\lam)-p(\lam_0)|
\le 0.5 |\lam|^{n+1}.
$$
Then for all $\lam\in B_\delta(\lam_0)$ we have
$p(\lam) \in \lam^{n+1} F$, hence $\lam\in \Mk$ by Lemma~\ref{lem-inclu}. 
\QED

\medskip

\noindent {\bf Example.}
It is not difficult to find specific
elements of $\Mk_0$ in $\Int(H)$. Using Bandt's algorithm for drawing
the set $\Mk$, for any ``black'' point $\lam_1$
in Figure 1, we can find a polynomial $p$ with coefficients in $\{0,\pm 1\}$
such that $|p(\lam_1)|$ is small. Then $p$ is likely to have a zero near
$\lam_1$. 

For example, we may take $\lam_1 = .14 + .68 i \in \Int(H)$ and the
polynomial $p(z)=1 + z^2 + z^3 - z^4 - z^5 + z^7$ (found with the help of
Bandt's algorithm), which has a zero
$\lam_0 \approx 0.141964+ 0.677696i$.
From the proof of Theorem~\ref{th-main} it follows that if $\lam\in \Int(H)$ and
$|p(\lam)-p(\lam_0)| \le 0.5 |\lam|^8$, then  $\lam\in \Mk$. We can
estimate, assuming that $|\lam|\le 2^{-1/2}$:
\be
|p(\lam)-p(\lam_0)| & \le & \sum_{k=1}^7 |\lam^k - \lam_0^k| \nonumber \\
& = & |\lam-\lam_0| \sum_{k=1}^7 |\lam^{k-1} + \lam\cdot \lam_0^{k-2} + 
\cdots + \lam_0^{k-1}| \nonumber \\
& \le & |\lam-\lam_0| \sum_{k=1}^\infty k \cdot 2^{-(k-1)/2} \nonumber \\
& = & |\lam-\lam_0|\cdot(1-2^{-1/2})^{-2}. \nonumber
\ee
It follows that the disc of radius $2\cdot 10^{-3}$ centered at 
$\lam_0$ lies in $\Mk$. (Indeed, $|\lam_0| > .692$, so $|\lam|>.69$ in this
disc, and $.69^8(1-2^{-1/2})^2/2>2\cdot 10^{-3}$.)

\medskip

\noindent{\em Proof of Lemma~\ref{lem-comput}.}
Let $R:=R_{a,b}$ be the rectangle with vertices $\pm a \pm ib$. We 
assume right away that $a\ge b$ and $a>1$.
Let $\lam = \xi + i\eta$. We assume that $\xi,\eta>0$, $|\lam|^2 \ge 
\frac{1}{3}$ and
$\xi \le |\lam|^2$ (this is certainly true in the set $H$). The
condition $R\subset \lam R \cup (\lam R -1) \cup (\lam R + 1)$ is
equivalent to 
\be \label{eq-cover}
\lam^{-1}R \subset R \cup (R-\lam^{-1}) \cup
(R + \lam^{-1}).
\ee
Figure 5 will help us write down sufficient conditions for (\ref{eq-cover});
it shows the case of $\lam = 0.1 + 0.68i$, $a = 1.35$, and $b=0.78$.
\begin{figure}[ht]
\epsfig{figure=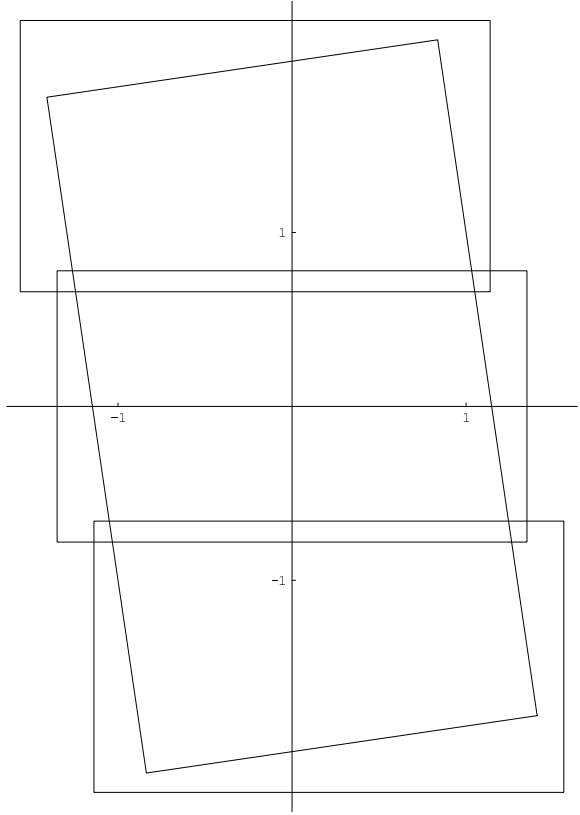,height=12cm}
\caption{$\lam^{-1}R_{a,b} \subset R_{a,b} \cup (R_{a,b}-\lam^{-1}) \cup
(R_{a,b} + \lam^{-1})$}
\end{figure}
We claim that the following conditions imply (\ref{eq-cover}):

(C1) {\sl Overlapping:} for $R$ to overlap $R\pm \frac{1}{\lam}$ we need that
$
|\IM(1/\lam)| < 2b,
$
or equivalently,
\be \label{c1}
|\lam|^2 \ge \frac{\eta}{2b}\,.
\ee

(C2) {\sl Covering the short sides:} we want the vertices 
$\frac{a\pm ib}{\xi + i\eta}$ of $\frac{1}{\lam}R$
to lie in $R+\frac{1}{\lam}$; then by
symmetry, the vertices $\frac{-a\pm ib}{\xi + i\eta}$ will be in
$R-\frac{1}{\lam}$. The conditions are:
\be \label{c21}
-a + \RE(1/\lam) \le \RE[(a\pm ib)/\lam] \le a + \RE(1/\lam),
\ee
\be \label{c22}
-b + \IM(1/\lam) \le \IM[(a\pm ib)/\lam] \le b + \IM(1/\lam).
\ee

(C3) {\sl Covering the long sides:} we want the upper intersection point
of $\partial (\lam R)$ and $\partial(\lam R+1)$ to be above or on 
the line $y=b$.
Then by symmetry the same is true for the intersection point
of $\partial (\lam R)$ and $\partial(\lam R-1)$, and the upper side of
$R$ is covered. Again by symmetry, the lower side will be covered as well.
With the help of Figure 5 we can see that this intersection point is
$\lam(a + i(b+\IM(1/\lam))$, so we get the condition
\be\label{c3}
\IM(\lam(a + i(b+\IM(1/\lam))) \ge b.
\ee

After some algebraic manipulations we get the following four inequalities
(it turns out that out of the eight inequalities in (\ref{c21}), (\ref{c22}) 
four are
always true by our assumptions, and two more follow from the remaining two).

(A1) $|\lam|^2 \ge \frac{\eta}{2b}$;

(A2) $a\eta + b\xi \ge b + \frac{\xi\eta}{|\lam|^2}$;

(A3) $a\xi + b\eta \le a|\lam|^2 + \xi$;

(A4) $a\eta + b\xi \le b|\lam|^2 + \eta$.

We want to find for which $(\xi,\eta)$ there exist $(a,b)$ such 
that (A1)--(A4) are satisfied. 
It is natural to look for the solution by setting some inequalities to
be equalities. 

First, set (A2) and (A3) to be equalities and solve for $a$ and $b$ to obtain
\be \label{eq-ab1}
a = 1 + \frac{(|\lam|^2-\xi)^2}{|\lam|^2(1-|\lam|^2)},\ \ \ \ \
b = \frac{\eta (|\lam|^2-\xi)}{|\lam|^2(1-|\lam|^2)}\,.
\ee
The inequalities (A1) and (A4) reduce to
\be \label{eq-xi1}
|\lam|^2 - \frac{|\lam|\sqrt{1-|\lam|^2}}{\sqrt{2}} 
\le \xi \le \frac{3|\lam|^2-1}{2}\,.
\ee

Next, set (A3) and (A4) to be equalities and solve for $a$ and $b$ to obtain
\be \label{eq-ab2}
a = \frac{\eta^2 + \xi(|\lam|^2-\xi)}{\eta^2 - (|\lam|^2-\xi)^2},\ \ \ \ \
b = \frac{\eta|\lam|^2}{\eta^2 - (|\lam|^2-\xi)^2}\,.
\ee
The inequality (A2) reduces to
\be \label{eq-xi2}
\xi \le |\lam|^2 - \frac{|\lam|\sqrt{1-|\lam|^2}}{\sqrt{2}}\,,
\ee
and (A1) follows from (A2) under our assumptions, whenever $|\lam|^2 \ge
\frac{1}{3}$. Combining (\ref{eq-xi1}), (\ref{eq-ab1}), (\ref{eq-xi2}),
and (\ref{eq-ab2}) yields that the desired $a,b$ may be found for
all $\lam = \xi + i\eta$ such that $\frac{1}{3} \le |\lam|^2 \le \half$ and
$0 \le \xi \le \frac{3|\lam|^2 -1}{2}$. This is equivalent to
$\lam\in H$.
It remains to check the initial assumptions on $a$ and $b$. 
It is immediate from (\ref{eq-ab1}) and
(\ref{eq-ab2}) that $a\ge 1$. 
To estimate $b$, we note that every $\lam\in H$ satisfies $|\lam|^2\le \half$,
so $\xi \le \frac{3|\lam|^2-1}{2} \le \frac{1}{4}$. On the other hand,
$\xi^2+\eta^2 = |\lam|^2 \ge \frac{1}{3}$. Thus, by (A1),
$$
b \ge \frac{\eta}{2|\lam|^2} \ge \eta \ge \sqrt{\frac{1}{3}-\xi^2}
\ge \sqrt{\frac{1}{3}-\frac{1}{16}} = \sqrt{\frac{13}{48}} > \frac{1}{2}\,,
$$
as desired.  It remains to check that $a\ge b$.
By (A2), $\frac{a}{b} \ge \frac{1-\xi}{\eta}$. 
But $|\lam|^2 = (\xi^2+\eta^2)\le \half$, hence
$\xi+\eta\le 1$, and so $\frac{a}{b} \ge 1$.
\QED

\medskip

\noindent {\bf Remarks.} 1.
One might hope to extend this method to show
that the whole set $\Mk$ is the closure of its interior
(as conjectured by Bandt \cite{bandt}). To do this, one would need to
prove that the attractor $A_\lam\{-1,0,1\}$
contains some neighborhood of the origin, with a radius locally uniformly
bounded below.
It is expected that for $|\lam|>3^{-1/2}$ (with $\lam$ non-real)
the attractor $A_\lam\{-1,0,1\}$ has positive
measure and non-empty interior, at least, for a typical $\lam$, since then
the similarity dimension of the IFS is greater than 2. However, this
remains an open problem.

2. In Lemma~\ref{lem-attr1}, it is enough that $F$ is covered
not just by the images under one iteration of the
IFS, but by images under any number of iterations;
the conclusion remains the same. Moreover, instead of
rectangles, one can consider other simple shapes and attempt to establish the
covering property. A similar method was employed in
\cite{IJK}. In that paper the authors used a disc rather than
rectangle for the set $F$. They considered $\lam = e^{-2\pi i/5}/\sqrt{2}$ and
had to use 5 iterations of the IFS to achieve the needed covering property. This
involved drawing $(3^6-1)/2$ circles and heavy computer use.


\section{The case of $|\lam|\ge 1/\sqrt{2}$}

Here we prove Theorems \ref{thm-pisot} and \ref{th-garsia} and
Proposition~\ref{prop-elem}.

\medskip

\noindent {\em Proof of Theorem~\ref{thm-pisot}.}
This is a straightforward extension of \cite{erdos1},
but we provide the argument for completeness.
Consider the Fourier transform
$$
\nulahat(\xi)=\int_{\R^2} e^{it\cdot\xi}\,d\nula(t)
$$
where $\xi, t \in \R^2$ and $t\cdot\xi$ is the dot product.
Considering $\xi$ and $t$ as complex numbers we have $t\cdot\xi=
\RE(t\ov{\xi})$. Then, by the 
independence of the random series defining $\nula$,
\be \label{nulahat}
\nulahat(\xi) = \prod_{n=0}^\infty
\widehat{\shalf(\delta_{\lam^n} + \delta_{-\lam^n})}(\xi)
= \prod_{n=0}^\infty\cos(\RE(\lam^n \ov{\xi})).
\ee

Let $\theta$ be a complex Pisot number, with $|\theta|^2 \in (1,2)$.
Then $\theta$ has degree at least three, since $|\theta|^2 =\theta\ov{\theta}$
is not an integer, and the product of all roots of the minimal polynomial
is an integer. 
Let $\Th_3,\ldots,\Th_m$ be the algebraic
(Galois) conjugates of $\Th$ other than $\ov{\Th}$. We know that
$|\Th_j|< 1$ for $j\ge 3$.
Observe that $\Th$ is an algebraic unit, since
the constant term of the minimal polynomial $(-1)^m |\theta|^2 \prod_{j=3}^m 
\Th_j$ is less than 2 in modulus, so it must be equal to 1.
Since 
$\Th^n + \ov{\Th}^n + \sum_{j=3}^m \Th_j^n \in \Z$ for all $n\ge 1$, we have
\be \label{eq-pisot1}
\dist(2\RE(\Th^n),\Z) \le c \rho^n,\ \ n\ge 1,
\ee
for some $\rho \in (0,1)$ and $c>0$.
Recall that $\Th= 1/\lam$. We have
\be \label{eq-erd1}\ \ \ \ \ \ \ 
\nulahat(2\pi\ov{\Th}^N) = \prod_{n=-\infty}^N \cos(2\pi\RE(\Th^n))=
\prod_{n=1}^N \cos(2\pi \RE(\Th^n)) \cdot \nulahat(2\pi). 
\ee
We claim that $\nulahat(2\pi\ov{\Th}^N) \ne 0$ for $N\in \Nat$.
Indeed, since $\Th$ is an algebraic unit, $\Th^{-1}$ is also an
algebraic unit, and therefore, $\Th^n$ is an algebraic integer for all
$n\in \Z$. Now
$2\RE(\Th^n)=\Th^n + \ov{\Th}^n$ is also an algebraic integer
(as a sum of algebraic integers) for all $n\in \Z$, 
so it cannot be equal $k + \half$, for
some integer $k$.  It follows that $\cos(2\pi\RE(\Th^n))\ne 0$ for
all $n\in \Z$, which implies the claim.

Now (\ref{eq-erd1}) and (\ref{eq-pisot1}) imply
$$
|\nulahat(2\pi\ov{\Th}^N)| \ge \prod_{n=1}^\infty |\cos(c\rho^n)|
\cdot|\nulahat(2\pi)|
=:\delta>0,
$$
for all $N\ge 1$.
Thus, $\nulahat(\xi) \not\to 0$ as $|\xi|\to \infty$. By the Riemann-Lebesgue
Lemma, $\nula$ is not absolutely continuous, so by the Law of Pure Type,
$\nula$ is singular. \qed

\medskip

{\em Proof of Theorem~\ref{th-garsia}.}
This is a straightforward extension of \cite{garsia1},
but we provide the argument for completeness.
For $n \ge 1$ let
$$
\Ak_\lam^{(n)}:= \left\{ \sum_{k=0}^{n-1} a_k \lam^k:\ a_k = \pm 1 \right\}
$$
We will prove the following two claims:
\be \label{gars2}
\# \Ak_\lam^{(n)} = 2^n.
\ee
\be \label{gars3}
\exists\,c>0,\ \forall\,n\ge 1,\ |x-y| \ge c\cdot 2^{-n/2}\ \ \ \mbox{for all}\
\ x,y\in \Ak_\lam^{(n)},\ x\ne y.
\ee
First we deduce the
desired statement follows from (\ref{gars2}) and (\ref{gars3}).
Absolute continuity of $\nula$ with a bounded density will follow
if we can show that 
\be \label{hold}
\nula(B_r(z))\le \const\cdot r^2\ \ \ \mbox{for all}\ z\in \Compl
\ee
for all $r>0$.
Let $C_1>0$ be such that $A_\lam \subset B_{C_1}$.
Clearly, it is enough to establish (\ref{hold}) 
for $r=C_1|\lam|^n$, for all $n\ge 1$. Fix $n\ge 1$.
By self-similarity, for any Borel set $E$,
$$
\nula(E) = \sum_{x\in \Ak_\lam^{(n)}} \nula\left(\frac{E-x}{\lam^n}\right)
\cdot 2^{-n}.
$$
Thus, the measure
$\nula$ is a sum of $2^n$ ``pieces,'' each having the measure $2^{-n}$.
These pieces are supported on $\lam^n A_\lam + x$, for $x\in \Ak_\lam^{(n)}$.
Note that the supports of the pieces lie 
in $B_{C_1|\lam|^n}(x)$. It follows that
for any $z\in \Compl$ we have 
$\nula(B_{C_1|\lam|^n}(z)) \le N\cdot 2^{-n}$, where $N$ is the 
number of points in $\Ak_\lam^{(n)}$ that lie in $B_{2C_1|\lam|^n}(z)$.
Since the separation
between those points is at least $c\cdot 2^{-n/2}$,
we have $N (c^2/4) \cdot 2^{-n} \le 4C_1^2|\lam|^{2n}$, whence
$\nula(B_r(z)) \le (16/c^2)r^2$, as desired. \qed

\medskip

{\em Proof of} (\ref{gars2}).
Suppose that there is a point in $\Ak_\lam^{(n)}$ having two different
representations. Then there is a non-trivial
polynomial $p(z) = a_0 + a_1 z + \cdots
+ a_{n-1} z^{n-1}$,
with $a_k\in \{-1,0,1\}$, such that $p(\lam)=0$.
The polynomial $q(z) = z^{n-1}p(z^{-1})$ vanishes at $\Th= \lam^{-1}$, hence
the minimal polynomial of $\Th$ divides $q$. This is a contradiction, since
$q$ has all coefficients of modulus less than or equal to one and 
the minimal polynomial of $\Th$ 
has the constant term $\pm 2$. \qed

\medskip

{\em Proof of} (\ref{gars3}).
Let $x,y\in \Ak_\lam^{(n)}$, $x\ne y$. We have
\be \label{gars4}
\shalf (x-y) = a_0 + a_1 \lam + \cdots + a_{n-1} \lam^{n-1} =: p(\lam),
\ee
where $a_i\in \{-1,0,1\}$ are not all zeros. Let $q(z)= z^{n-1}p(z^{-1})$.
Let $\Th_1=\Th$, $\Th_2=\ov{\Th}$, and let $\Th_3,\ldots,\Th_m$ be
the remaining algebraic conjugates of $\Th$.
The product $\prod_{j=1}^m q(\Th_j)$
is an integer, since it is a value of a symmetric polynomial on the roots
of the minimal polynomial of $\Th$.
We know that $q(\Th)\ne 0$
by the proof of (\ref{gars2}), so $q(\Th_j)\ne 0$ for all $j\le m$, since
$q$ has integer coefficients. Thus,
\be \label{gars5}
\left| \prod_{j=1}^m q(\Th_j) \right| \ge 1.
\ee
Since $|\Th_j|>1$ by assumption, we have
$$
|q(\Th_j)| \le |\Th_j|^{n-1} + |\Th_j|^{n-2} + \cdots <
\frac{|\Th_j|^n}{|\Th_j|-1}\,.
$$
The constant term of the minimal
polynomial is $\pm 2$, so
\be \label{gars1}
\left|\prod_{j=3}^m \Th_j \right| = 2 |\Th|^{-2}.
\ee
In view of (\ref{gars5}) and (\ref{gars1}),
\be
|q(\Th)|^2 & = & q(\Th) q(\ov{\Th}) \ge \prod_{j=3}^m |q(\Th_j)|^{-1}
\nonumber \\
           & \ge & \prod_{j=3}^m (|\Th_j|-1) |\Th_j|^{-n} \nonumber \\
           & =   & 2^{-n} |\Th|^{2n} \prod_{j=3}^m (|\Th_j|-1). \nonumber
\ee
Now (\ref{gars4}) implies
$$
|x-y| = 2|p(\lam)| = 2 |\Th|^{-n+1} |q(\Th)| \ge c\cdot 2^{-n/2},
$$
as desired, where $c= 2|\Th|(\prod_{j=3}^m (|\Th_j|-1)^{1/2}$.
\qed

\medskip

\noindent {\em Proof of Proposition~\ref{prop-elem}(i).}
This statement is immediate from the following lemma.
Recall that $R_{a,b}$ is the rectangle with the vertices $\pm a \pm ib$.

\begin{lemma} \label{lem-cover2}
For any $\lam \in \Omega$
there exist $a>1, b>2^{-1/2}$ such that
$R_{a,b} \subset (\lam R_{a,b} -1) \cup
(\lam R_{a,b} +1)$.
\end{lemma}

The proof of this lemma is analogous to that of 
Lemma~\ref{lem-comput}, so we omit it. 

\medskip

\noindent {\em Proof of Proposition~\ref{prop-elem}(ii)}
is the same as in \cite[p.\,424]{IJK}; we include it for completeness.
We prove the first statement; the second one is proved similarly.
Let $\lam = re^{\pi i \frac{m}{n}}$,
with  $GCD(m,n)=1$. By (\ref{eq-ss2}), 
\be \label{eq-elem1}
A_\lam = A_{\lam^n} + \lam A_{\lam^n} + \cdots + \lam^{n-1} A_{\lam^n}.
\ee
Note that $\lam^n = \pm r^n$ and $A_{z} = A_{-z}$, so $A_{\lam^n} = 
A_{r^n}\subset \R$. Moreover,
for $r \in [2^{-1/n},1)$ we have $r^n \ge \half$, so $A_{r^n}$ is a line
segment of length $a = \frac{2}{1-r^n}$. Now we see that (\ref{eq-elem1})
represents a Minkowski sum of $n$ line segments making the angles
$\frac{\pi k}{n}$, $k=0,\ldots,n-1$, with the horizontal. This implies the
desired statement.
\qed


\section{Results for a typical $\lam$}

\noindent {\em On the proof of Proposition~\ref{prop-trans}.} 
In \cite[Thm 2]{bbbp} it is proved that
if a power series $1 + \sum_{n=1}^\infty a_n z^n$, with $a_n \in [-1,1]$,
has $k$ roots (counting with multiplicities)
in a disc of radius $r$, then
\be \label{eq-trans}
r \ge k^{-1/2k} \left(1 + \frac{1}{k}\right)^{-\half(1+1/k)}.
\ee
The presence of a
non-real double root $\lam$ implies that there are at least four roots
in the disc of radius $|\lam|+\eps$, since the power series has real
coefficients and $\overline{\lam}$ has to be a double root as well.
Substituting $k=4$ into the formula yields the desired estimate. \qed

\medskip

\noindent {\bf Remark.} The proof of (\ref{eq-trans}) is remarkably simple;
it is based on  Jensen's Formula, concavity of the logarithm, and 
Parseval's Formula. The bound 
$2\cdot 5^{-5/8} \approx .73143$ obtained this way is quite good: 
Pinner \cite{pinner}
gave an example of a power series with coefficients in $[-1,1]$ having
a double zero of modulus $\approx .75361$.

\medskip

Using Proposition~\ref{prop-trans}, all the theorems from Subsection 2.1.4
follow by standard methods. We only indicate 
how to deduce Theorem~\ref{th-meas}(ii).

\medskip

By (\ref{nulahat}), for any $k \ge 2$,
\be \label{nuhat}
\nulahat(\xi) = \nuhat_{\lam^k} (\xi) \nuhat_{\lam^k} (\lam \xi) \cdot
\ldots \cdot \nuhat_{\lam^k} (\lam^{k-1} \xi).
\ee
Suppose that $\nu_{\lam^k}$ has a density in $L^2(\R^2)$. Then $\nuhat_{\lam^k}
\in L^2(\R^2)$ by Plancherel's Theorem, 
hence $\nulahat \in L^{2/k}(\R^2)$ by (\ref{nuhat}).
Since $\nula$ is a probability measure, $\nulahat$ is also bounded, so
$\nulahat \in L^1(\R^2)$. Applying the Inverse Fourier Transform, we
conclude that $\nula$ is absolutely continuous with a continuous density. 
Thus, for any $k\ge 2$,
$$
\frac{d \nu_{\lam^k}}{dx} \in L^2(\R^2) \ \Rightarrow \ 
\frac{d\nula}{dx} \in C(\R^2).
$$
Now Theorem~\ref{th-meas}(i) implies Theorem~\ref{th-meas}(ii). \qed

\medskip

\subsection*{Acknowledgement}
We are indebted to Christoph Bandt for motivating questions, many useful
discussions, and help with the figures. Thanks also to David Boyd and
David Garth for valuable information concerning complex Pisot numbers.


\bibliographystyle{amsplain}

\end{document}